\documentclass[11pt]{article}

\usepackage[top=4cm, bottom=2cm, left=2.2cm, right=2.2cm]{geometry}
\usepackage{subfig}  
\usepackage{amssymb}
\usepackage{graphicx}
\usepackage{amsmath}
\usepackage[titletoc]{appendix}
\usepackage{appendix}
\usepackage{pifont}
\usepackage{bbm}
\usepackage{amsfonts}
\usepackage{mathrsfs}
\usepackage{amssymb}
\usepackage{amsmath}
\usepackage{a4wide}
\usepackage{stmaryrd}
\usepackage{txfonts}
\usepackage{tabularx}
\usepackage{booktabs}
\usepackage{ctable}

\usepackage{latexsym}

\let\er\eqref

\newcommand{\R}{{\mathbb R}}

\newtheorem{thm}{Theorem}
\newtheorem{lemma}[thm]{Lemma}
\newtheorem{proposition}[thm]{Proposition}
\newtheorem{remark}[thm]{Remark}

\begin{document}
\fontfamily{ptm}\fontseries{sb}\selectfont
\title{
{\bf\Large Nonlocal nonlinear reaction preventing blow-up in supercritical case of chemotaxis system}
 \\
}
\author{ Shen Bian \footnote{Beijing University of Chemical Technology, 100029, Beijing. Universit\"at Mannheim, 68131, Mannheim. Email: \texttt{bianshen66@163.com}. Partially supported by National Science Foundation of China (Grant No. 11501025) and the Alexander von Humboldt Foundation.}
  \and Li Chen\footnote{Universit\"at Mannheim, 68131, Mannheim. Email: \texttt{chen@math.uni-mannheim.de}. Partially supported by DFG Project CH 955/3-1.}
  \and Evangelos A. Latos\footnote{Universit\"at Mannheim, 68131, Mannheim. Email: \texttt{evangelos.latos@math.uni-mannheim.de}. Partially supported by DFG Project CH 955/3-1.}
}
\date{}
\maketitle
\begin{center}
{\bf\small Abstract}

\vspace{3mm}\hspace{.05in}\parbox{5in} {
This paper is devoted to the analysis of non-negative solutions for the chemotaxis model with nonlocal nonlinear source in bounded domain. The qualitative behavior of solutions is determined by the nonlinearity from the aggregation and the reaction. When the growth factor is stronger than the dampening effect, with the help of the nonlocal nonlinear term in the reaction, for appropriately chosen exponents and arbitrary initial data, the model admits a classical solution which is uniformly bounded. Moreover, when the growth factor has the same order with the dampening effect, the nonlocal nonlinear exponents can prevent the chemotactic collapse.}
\end{center}

\noindent
{\it \footnotesize \textbf{Key words.}} {\footnotesize 
Chemotaxis model, Fisher-KPP model, Global existence,
Nonlocal reaction}

\section{Introduction}
\def\theequation{1.\arabic{equation}}\makeatother
\setcounter{equation}{0}
\def\thetheorem{1.\arabic{thm}}\makeatother
\setcounter{thm}{0}


The Keller-Segel model in Chemotaxis was originally introduced by Keller and Segel \cite{KS1,KS2} to describe the characteristic movement of cells, the cells can move toward the increasing signal concentration or can be repulsive by the signal concentration. From then on, mathematical models to describe chemotaxis have been widely proposed in the last few years. The simplest version contains the competition among the diffusion, reproduction and the nonlocal aggregation satisfying \cite{bp07}
\begin{align} \label{star11}
\left\{
      \begin{array}{ll}
       u_t=\Delta u-\chi \nabla \cdot (u^{\sigma} \nabla c)+f(u), &  x \in \Omega, t>0, \\
     \tau c_t-\Delta c+ c=u^{\xi}, & x \in \Omega, t>0, \\
     u(x,0)=u_0(x) \geq 0, &  x \in \Omega.
      \end{array}\right.
\end{align}
In the modelling, $\Omega$ is either a bounded domain in $\R^n$ or the whole space. In the context of biological aggregation, $u(x,t)$ represents the bacteria density, $c(x,t)$ is the chemical substance concentration. The reaction term describes the reproduction rate of the bacteria where the resources of the environment can be consumed either locally or nonlocally. When chemicals diffuse much faster than cells \cite{JL92}, \er{star11} can be reduced into parabolic-elliptic model, i.e.
\begin{align} \label{star00}
\left\{
      \begin{array}{ll}
       u_t=\Delta u-\chi \nabla \cdot (u^{\sigma} \nabla c)+f(u), &  x \in \Omega, t>0, \\
     -\Delta c+ c=u^{\xi}, & x \in \Omega, t>0.
      \end{array}\right.
\end{align}
In the following we will report some of the related previous results on \er{star00} in terms of $f(u)$.

\er{star00} with $f(u) \equiv 0$ is the classical Keller-Segel model which expresses the random movement(brownian motion) of the cells with a bias directed by the chemoattractant concentration \cite{bp07}. This system has been widely studied, such as \cite{HV96,JL92,Na95} and the references therein. It's proved that for the following problem with Neumann boundary condition
\begin{align*}
& u_t=\Delta u-\chi \nabla \cdot (u \nabla c),  \\
& -\Delta c=u-1,
\end{align*}
blow-up never occurs in one dimension \cite{Na95}. While in two dimensions, there exists a threshold number for the initial data that can separate global existence and finite time blow-up \cite{JL92}.

When $f(u)\neq 0$, the logistic growth including the consumption of resources around the environment is taken into account in chemotaxis models. There are quite a number of works handling such type of model with logistic growth describes the situation where the influence of nonlocal terms is omitted. Here we can only list some of the results which are closely related to our model.

For $\sigma=\xi=1,$ the authors in \cite{TW07} proved that model \er{star00} with
\begin{align}
f(u) \le a-b u^2,~~u \ge 0
\end{align}
possesses a global bounded classical solution for either $n \le 2$ or $n \ge 3$ and $b>\frac{n-2}{2}\chi$. In addition, for all $n \ge 1, b>0$
and arbitrary initial data there exists at least one global weak solution given by $f(u) \ge -c_0(u^2+1),~\forall u>0$ with some $c_0>0.$

For more general case, in \cite{WangMuChZ14} the authors considered the model
\begin{align*}
& u_t=\nabla \cdot (D(u) \nabla u)-\chi \nabla \cdot (u \nabla c)+f(u), \\
& -\Delta c+c=u.
\end{align*}
Here $f(u)$ is smooth satisfies $f(0) \ge 0$ and
$$f(u) \le a-b u^{\gamma}$$
for all $u \ge 0$ with $a \ge 0, b>0$ and $\gamma>1$. $D(u) \in C^2\left([0,\infty)\right)$ and there exist some constants $c_D>0$ and $m \ge 1$ such that
$D(u) \ge c_D u^{m-1}$ for all $u >0$ as well as $D(u)>0$ for all $u \ge 0$. They proved that if $\gamma \ge 2$ and $b>b^*$ where
$$
b^*=
\begin{cases}
	\frac{(2-m)n-2}{(2-m)n}\chi,\quad
&if \quad m<2-2/n,\\
0,\quad &if \quad m\geq2-2/n,
\end{cases}
$$
or $\gamma\in(1,2)$ and $m>2-2/n$, then the model has a unique nonnegative classical solution which is global and bounded.

In \cite{Galakhov:2016bka} the authors considered
\begin{align*}
 &u_t=\Delta u-\chi \nabla \cdot (u^{\sigma} \nabla c)+\mu u(1-u^\alpha), \\
 &-\Delta c+c=u^{\xi}
\end{align*}
with $\sigma \ge 1, \xi \ge 1.$ If $\alpha> \sigma+\xi-1$ or
\begin{align}
\alpha=\sigma+\xi-1~~\mbox{and}~~ \mu> \frac{n \alpha-2}{n \alpha+2(\sigma-1)}\chi,
\end{align}
then for sufficiently smooth initial data, there exists a unique global solution of the model. Afterwards in \cite{Hu:2017kd} it was proved that the same result for the above problem holds true even for the case
\begin{align}
\alpha=\sigma+\xi-1~~\mbox{and}~~ \mu= \frac{n \alpha-2}{n \alpha+2(\sigma-1)}\chi.
\end{align}

Logistic growth described by nonlocal terms has been investigated in recent years. For example, \cite{NegTel13} focused on the parabolic-elliptic system with linear competitive effect
\begin{align}\label{NT13}
\left\{
      \begin{array}{ll}
       u_t=\Delta u-\chi \nabla \cdot (u \nabla c)+u\left( a_0-a_1 u-\frac{a_2}{|\Omega|} \int_{\Omega} u dx \right), \\
     -\Delta c+ c=u+g
      \end{array}\right.
\end{align}
where $a_0,a_1>0,\chi>0,a_2 \in \R$ and $g$ is a uniformly bounded function. Actually, as the population grows, the competitive effect of local term $a_1 u$ is more influential than the nonlocal term $\int_{\Omega} u dx$. In this case, the reaction term $f(u)$ behaves like $u(a_0-a_1 u)$, following a comparison argument based on upper and lower solutions defined by a global solution of an ordinary system, the authors showed that if $a_1>2 \chi+|a_2|$, then $\|u-\frac{a_0}{a_1+a_2}\|_{L^\infty(\Omega)} \to 0$ as time goes to infinity.

Since there is a fertile area for this research, it's difficult to cover all the important results, we refer the interested readers to \cite{He:2016fz,Issa:2016uq,NakOsa11,NegTel13,SzyRoLaCha09,WL17,ZhLi15,Zh15}.

In this paper, our main purpose is concerned with the effect of the nonlocal nonlinear source for the solutions of non-degenerate model \er{star00}. Before our main result, let's present some interesting aspects connected to a deeper understanding of our paper.

Firstly, in the nonlocal term of \er{NT13}, as the population grows, the competitive effect of the local term becomes more influential than the nonlocal term, and the effect of the total mass can be ignored compared with the local term, hence the influence of nonlocal term is still unknown. As it was stated in \cite{NegTel13} that ``it seems to conjecture that the dampening effect of the nonlocal terms might lead to an even more effective homogenization, this case provides no information about the asymptotic behavior by a comparison method." In addition, logistic growth described by nonlocal terms has been used in a competitive system modelling cancer cells behavior which considers the influence of the surrounding area of a cell to replicate itself \cite{NegTel13,SzyRoLaCha09} and it can also describe Darwinian evolution of a structured population \cite{KPP} or nuclear reaction process \cite{HY95,WW96}. Therefore the effect of the nonlocal term on the diffusion-aggregation-reaction equation is also very attractive.

Secondly, in light of the known research, the available analytical results on chemotaxis with logistic sources mostly concentrate on local reaction term, i.e.
\begin{align}\label{presum}
u_t=\Delta u-\chi \nabla u^{\sigma} \cdot \nabla c-\chi u^\sigma c+ \chi u^{\sigma+\xi}+a u-b u^\alpha,
\end{align}
model \er{presum} possesses a global classical solution with the fact that either the dampening effect is stronger than the growth factor \cite{Galakhov:2016bka,WangMuChZ14}, i.e. $\alpha>\sigma+\xi$, or the dampening effect has the same order with the growth factor \cite{NegTel13,TW07}, i.e. $\alpha=\sigma+\xi$, combining some constraints on the coefficients $b$ and $\chi$. To the best of our knowledge, when the growth factor is stronger than the dampening effect, whether the non-degenerate model \er{star00} admits a global solution is still open.

Therefore, in this paper, in order to detect the influence of nonlocal term on the behavior of solutions, without loss of generality, the coefficients of the dampening term and the growth factor are fixed to be $b=1, \chi=1$ and $\sigma, \xi$ in model \er{presum} are constrained to be $\sigma=1, \xi=1$. More precisely, we will study the following chemotaxis system with nonlocal reaction
\begin{align} \label{nkpp00}
\left\{
      \begin{array}{ll}
      \ u_t=\Delta u-\nabla \cdot (u \nabla c)+f(u), &  x \in \Omega, t>0, \\
      -\Delta c+c=u, & x \in \Omega, t>0, \\
    \nabla u \cdot \nu=\nabla c \cdot \nu=0,  & x \in \partial \Omega, t>0, \\
     u(x,0)=u_0(x)\geq 0, &  x \in \Omega,
      \end{array}\right.
\end{align}
where $\Omega$ is a smooth bounded domain in $\R^n(n \ge 3)$ and $\nu$ is the outer unit normal vector on $\partial \Omega$, the reaction term is taken to be
$$
f(u)=u^\alpha \left(1-\int_{\Omega} u^\beta dx\right)
$$
with $\alpha \ge 1, \beta>1$. The initial data is assumed to be
\begin{align}\label{initial0}
u_0 \ge 0, ~u_0 \in C^{\theta}(\overline{\Omega})\mbox{~for some ~} \theta \in (0,1).
\end{align}
Actually, model \er{nkpp00} can be rewritten as
\begin{align*}
u_t=\Delta u-\nabla u \cdot \nabla c-u c+ u^{2}+u^\alpha-u^\alpha \int_{\Omega} u^\beta dx.
\end{align*}
If the dominated growth factor $u^2$ is stronger than the dampening effect $u^\alpha$, the nonlocal term $\int_{\Omega} u^\beta dx$ can help preventing the chemotactic collapse. Precisely, our result is the following:
\begin{thm}\label{thm0}
Let $n \ge 3, \alpha \ge 1, \beta>1$, $u_0$ satisfies \er{initial0}.
If
\begin{align}\label{con01}
2 \le \alpha <1+2 \beta/n
\end{align}
or
\begin{align}\label{con02}
\alpha<2 \mbox{~~~and~~~} \frac{n+2}{n} \left(2-\alpha\right) < 1+2\beta/n-\alpha,
\end{align}
then the problem \er{nkpp00} admits a unique global classical solution which is uniformly bounded. Besides, the following estimate for any $t>0$ holds true,
\begin{align}
\|u(\cdot, t)\|_{L^\infty(\Omega)} \le C(\|u_0\|_{L^1(\Omega)},\|u_0\|_{L^\infty(\Omega)}).
\end{align}
\end{thm}

\begin{remark}
It is well known, for example in \cite{JL92}, that the solution of chemotaxis system without reaction blows up in finite time for large initial data. Comparing this result with Theorem \ref{thm0} shows that an appropriate nonlocal nonlinear dampening effect could give a global in time solution without any restriction on the initial data. The condition of Theorem \ref{thm0} implies $\beta>n/2$ can prevent chemotactic collapse. However, whether there exists a blow-up solution to model \er{nkpp00} under the assumption that $\beta \le n/2$ for higher dimension is still open.
\end{remark}

\begin{remark}
For $\beta \le n/2$ and the production term $u$ in the second equation of \er{nkpp00} is replaced by sub-linear term $u^\xi$ with certain $\xi<1$, if the assumption of Theorem \ref{thm0} is replaced by $1+\xi<1+2 \beta/n$, then model \er{nkpp00} asserts boundedness of solutions.
\end{remark}


In this paper, Section \ref{sec2} is devoted to the global solutions of model \er{nkpp00}, with that target the proof of the local existence and the key a priori estimates are presented. Precisely, some preliminary inequalities which are important for our proof are given in subsection \ref{subsec0}. Subsection \ref{subsec1} applies the Schauder fixed point theorm to show the local existence of classical solutions and blow-up criterion, where a careful application of Maximum principle is used in building up the mapping. Furthermore, $L^k(1 \le k<\infty)$ estimates are obtained by applying Sobolev type of inequalities, where $-u^\alpha \int_{\Omega} u dx$ from the reaction term plays a key rule so as to control the aggregation and nonlinear growth.



\section{Global bounded solution}\label{sec2}
\def\theequation{2.\arabic{equation}}\makeatother
\setcounter{equation}{0}
\def\thetheorem{2.\arabic{thm}}\makeatother
\setcounter{thm}{0}

This section is devoted to prove the global existence of solutions. Throughout the proof, we use the following exponent arising from Sobolev inequality \cite{lieb202}
\begin{align}\label{p}
      p=\frac{2n}{n-2}, ~~n \ge 3
\end{align}
for the convenience of calculations. Without loss of generality, we suppose $|\Omega|=1$.

\subsection{Preliminary}\label{subsec0}
Before showing the global existence, we need the following preparations. These lemmas have been proved in \cite{BL16,BL14}.
\begin{lemma}(\cite{BL16}) \label{lemmainterpolation}
Let $p$ is expressed by \er{p}, $1 \le r<q<p$ and $\frac{q}{r}<\frac{2}{r}+1-\frac{2}{p}$, then for $v\in H^1(\Omega)$ and $v \in L^r(\Omega)$, it holds
\begin{align}
\|v\|_{L^q(\Omega)}^q \le C(n) \left(C_0^{-\frac{ \lambda q}{2-\lambda q}}  + C_1^{-\frac{ \lambda q}{2-\lambda q}} \right)\|v\|_{L^r(\Omega)}^{\gamma} + C_0 \|\nabla v\|_{L^2(\Omega)}^2+ C_1 \|v\|_{L^2(\Omega)}^2,~~n \ge 3. \label{inter}
\end{align}
Here $C(n)$ are constants depending on $n$, $C_0,C_1$ are arbitrarily positive constants and
\begin{align}
\lambda=\frac{\frac{1}{r}-\frac{1}{q}}{\frac{1}{r}-\frac{1}{p}} \in (0,1),~~
\gamma=\frac{2(1-\lambda) q}{2-\lambda q}=\frac{2\left(1-\frac{q}{p}\right)}{\frac{2-q}{r}-\frac{2}{p}+1}.
\end{align}
\end{lemma}
For the $L^\infty$ estimates, we need the following inequality
\begin{lemma}(\cite{BL14})\label{odeinfinity}
Assume $y_k(t) \ge 0,~k=0,1,2,...$ are $C^1$ functions for $t>0$
satisfying
\begin{align}\label{yk}
y_k'(t) \le -y_k+a_k \left( y_{k-1}^{\gamma_{1}}(t)
+y_{k-1}^{\gamma_{2}} (t) \right),
\end{align}
where $a_k=\bar{a} b^{r_0 k}>1$ with $\bar{a},r_0,b$ are positive bounded
constants and $0<\gamma_{2}<\gamma_{1} \le b$. Assume also that
there exists a bounded constant $K \ge 1$ such that $y_k(0) \le
K^{b^k}$, then
\begin{align}
y_k(t) \le (2\bar{a})^{\frac{b^k-1}{b-1}} b^{r_0 \left(
\frac{b(b^k-1)}{(b-1)^2}-\frac{k}{b-1} \right)} \max \left\{
\sup_{t \ge 0} y_0^{b^k}(t), K^{b^k} \right\}.
\end{align}
\end{lemma}

\begin{remark}
Lemma \ref{odeinfinity} is an application of the Ghidalia's lemma (see \cite[Lemma 5.1]{T}).
\end{remark}
\subsection{Proof of Theorem \ref{thm0}} \label{subsec1}
We are now in a position to begin the study. In order to prove Theorem \ref{thm0}, we split the proof into three parts. Firstly in Proposition \ref{pro0}, we consider the local existence and uniqueness as well as the blow-up criterion of the classical solution. Then Proposition \ref{pro1} presents the a prior estimates which assure the uniformly boundedness of solutions. Finally we can directly obtain the global existence of the unique classical solution to close the proof of Theorem \ref{thm0}.

We firstly claim the result about the local existence of the classical solution to \er{nkpp00}.
\begin{proposition}\label{pro0}
Let $\alpha \ge 1.$ Assume $u_0\in C^{\theta}(\overline{\Omega})$ for some $\theta \in (0,1)$, then there exists a maximal existence time $T_{\max} \in (0,\infty]$ and a unique classical solution $u(x,t)$ to model \er{nkpp00} such that
\begin{align}
u \in C^{2+\delta,~1+\frac{\delta}{2}}\left( \overline{\Omega}\times(0,T_{\max}) \right),
\end{align}
where $\delta \in (0,1)$. Besides, if $T_{\max}<\infty,$ then
\begin{align}\label{blowup}
\displaystyle \lim_{t \to T_{\max}} \|u(t)\|_{L^\infty(\Omega)}=\infty.
\end{align}
\end{proposition}
\noindent\textbf{Proof of Proposition \ref{pro0}.} The local existence in time and blow-up criterion can be derived through a very standard demonstration. Here we refine the proof in spirit of \cite{TW07}. Firstly we construct a nonlinear ODE
\begin{align}
&\overline{u}_t=\overline{u}^{2}+\overline{u}^\alpha ,~~t \in (0,T) \label{ubar} \\
&\overline{u}(0)=\|u_0\|_{L^\infty(\Omega)}. \nonumber
\end{align}
Here $T=\frac{\tilde{T}_{\max}}{2}$ and $\tilde{T}_{\max}$ is the maximum existence time of $\overline{u}.$ Since $T<\tilde{T}_{\max},$ one has that $\overline{u} \ge 0$ is bounded in $[0,T]$. Here we denote $\displaystyle \max_{t \in [0,T]} \overline{u}(t)=L_0.$

Define the closed bounded subset
\begin{align}
S:=\big\{ \tilde{u} \in C_{x,t}^{\theta,\frac{\theta}{2}}(\overline{\Omega}\times [0,T])\big|~0 \le \tilde{u}(x,t) \le \overline{u}(t)\le L_0 \mbox{~in}~\Omega \times [0,T] \big\}.
\end{align}
We introduce a mapping $\Phi:~S \to S$ such that $\Phi(\tilde{u})=U.$ Here $U$ can be obtained from the following steps.

Firstly, we consider
\begin{align} \label{V00}
\left\{
      \begin{array}{ll}
  -\Delta V+V=\tilde{u}. ~&x \in \Omega,~0<t<T, \\
  \nabla V \cdot \nu=0, ~ &x \in \partial \Omega, ~0<t<T.
      \end{array}\right.
\end{align}
Since $\tilde{u}\in C^{\theta,\theta/2}(\overline{\Omega}\times [0,T]),$ by the theory of classical solutions to elliptic equations \cite[Theorem 8.34]{Trud} one can obtain that there is a unique solution
\begin{align}\label{V}
V(x,t) \in C^{2+\widehat{\theta},\frac{\widehat{\theta}}{2}}(\overline{\Omega}\times [0,T]),
\end{align}
then
\begin{align}\label{gradV}
\nabla V \in C^{1+\widehat{\theta},\frac{\widehat{\theta}}{2}}(\overline{\Omega}\times [0,T]).
\end{align}
In addition, by the maximum principle \cite[Theorem I 2.1]{LSU} one has
\begin{align}
0 \le V(x,t) \le \overline{u} \le L_0.
\end{align}
Secondly, we construct
\begin{align} \label{U00}
\left\{
      \begin{array}{ll}
      U_t=\Delta U- \nabla U \cdot \nabla V- U \Delta V +U \tilde{u}^{\alpha-1} \left( 1- \int_{\Omega} \tilde{u}^\beta dx \right), ~x \in \Omega,~0<t<T, \\
  \nabla U \cdot \nu=0, ~ x \in \partial \Omega, ~0<t<T, \\
  U \big|_{t=0}=u_0(x) \in C^{2+\delta}(\overline{\Omega}) \ge 0.
      \end{array}\right.
\end{align}
Since $V$ satisfies \er{V} and \er{gradV}, then there exists $\delta \in (0,1)$ depending on $\alpha,\beta,\theta$ such that the terms $\nabla V, \tilde{u}-V, \tilde{u}^{\alpha-1}( 1-\small{\int}_{\Omega} \tilde{u}^\beta dx )$ belong to $C^{\delta,\delta/2}\left( \overline{\Omega} \times (0,T] \right)$. Therefore,
from \cite[Theorem IV 5.4]{LSU} we have that \er{U00} has a unique classical solution $U(x,t) \in C^{2+\delta,1+\delta/2}(\overline{\Omega} \times [0,T])$ with $u_0 \in C^{2+\delta}(\overline{\Omega})$.

Now we will prove that $0 \le U \le \overline{u}(t)~~\mbox{in}~\overline{\Omega} \times (0,T].$ By the maximum principle \cite[Theorem I 2.1]{LSU} we obtain
$$
U \ge 0.
$$
On the other hand, $\overline{u}(t)$ satisfies
\begin{align*}
& \overline{u}_t-\Delta \overline{u}+\nabla\cdot( \overline{u} \nabla V)- \overline{u} \tilde{u}^{\alpha-1} \left(1-\int_{\Omega} \tilde{u}^\beta dx \right) \\
=&~\overline{u}_t+\overline{u} (V-\tilde{u})- \overline{u} \tilde{u}^{\alpha-1} \left(1-\int_{\Omega} \tilde{u}^\beta dx \right)    \\
\ge & ~ \overline{u}_t-\overline{u}^{2}-\overline{u}^{\alpha}=0,
\end{align*}
thus again applying the maximum principle with $U(0) \le \overline{u}(0)$ we have
\begin{align}
U(x,t) \le \overline{u}(t) \le L_0 \mbox{~in~} \overline{\Omega} \times [0,T].
\end{align}
Therefore $\Phi:~S \mapsto S$ is well defined. Since $U(x,t) \in C^{2+\delta,1+\delta/2}(\overline{\Omega}\times [0,T])$ is compactly embedding into $C^{\theta,\theta/2}(\overline{\Omega}\times [0,T]),$ thus $\Phi(S) \subset S$ is a relatively compact subset. Applying the Schauder fixed point theorem we have that there exists a fixed point of $\Phi$ which is the classical solution of \er{nkpp00}.

Let's mention that if $u_0 \in C^{\theta}(\overline{\Omega})$ for some $\theta \in (0,1),$ we can construct a sequence in $C^{2+\delta}(\overline{\Omega})$ that converges to $u_0$ in $C^{\theta}(\overline{\Omega})$. By inner $C^{2+\delta,1+\delta/2}$ regularity and compactness arguments, we can easily show that $u(x,t) \in C^{2+\delta,1+\delta/2}(\overline{\Omega}\times (0,T]).$

In addition, assume $u_1,u_2$ are two solutions of \er{nkpp00} in $(0,T]$, the multiplication \er{nkpp00}$_{u=u_1}-$\er{nkpp00}$_{u=u_2}$ by $u_1-u_2$ and the integration over $\Omega$ give that $\|u_1-u_2\|_{L^2(\Omega)} \equiv 0$ which assures the uniqueness of solutions.

Finally, parabolic regularity theory \cite[Theorem V 6.1]{LSU} follows that if $u$ is H\"{o}lder continuous, then the solution can be extended to the interval $[0,T_{\max})$ with $T_{\max} \le \infty$ and if $T_{\max}< \infty,$ \er{blowup} holds true by the standard arguments in \cite{HW2005}. Thus completes the proof. $\Box$ \\

The most important part to show the global existence is the following a priori estimates.

\begin{proposition}\label{pro1}
Let $\alpha \ge 1$, $u_0$ satisfies \er{initial0}. Let $u$ be any nonnegative classical solution of problem \er{nkpp00} within $0<t<T_{\max}.$ If either
\begin{align}\label{con01}
2 \le \alpha <1+2 \beta/n
\end{align}
or
\begin{align}\label{con02}
\alpha<2 \mbox{~~~and~~~} \frac{n+2}{n}(2-\alpha)< 1+2\beta/n-\alpha,
\end{align}
then the following estimate holds true that for any $0<t<T_{\max}$ and any $1 \le k < \infty$
\begin{align}
\int_{\Omega} u^k(\cdot,t) dx \le C\left( \|u_0\|_{L^1(\Omega)},\|u_0\|_{L^k(\Omega)}  \right).
\end{align}
Furthermore, the uniformly boundedness is obtained that for $0<t<T_{\max}$
\begin{align}\label{uniformbounded}
\|u(\cdot,t)\|_{L^\infty(\Omega)} \le C\left( \|u_0\|_{L^1(\Omega)},\|u_0\|_{L^\infty(\Omega)} \right).
\end{align}
\end{proposition}
\noindent\textbf{Proof of Proposition \ref{pro1}.} Beginning with a priori estimates we get the boundedness of $L^k$ norm for $1<k<\infty.$ Then using Lemma \ref{odeinfinity}, the uniformly boundedness of the solutions can be obtained by the iterative method.

\noindent {\bf Step 1 (A priori estimates).} It's obtained after multiplying \er{nkpp00} by $k u^{k-1}(k \ge 1)$ that
\begin{align}\label{star}
& \frac{d}{dt} \int_{\Omega}  u^k dx  + \frac{4(k-1)}{k}
\int_{\Omega} |\nabla u^{\frac{k}{2}} |^2 dx + k \int_{\Omega} u^\beta dx \int_{\Omega} u^{k+\alpha-1} dx  \nonumber\\
=& ~k \int_{\Omega} u^{k+\alpha-1} dx -k \int_{\Omega} u^{k-1} \nabla \cdot (u \nabla c) dx.
\end{align}
Recalling \er{nkpp00},
\begin{align}
-k \int_{\Omega} u^{k-1} \nabla \cdot (u \nabla c) dx= (k-1) \int_{\Omega} u^{k+1} dx - (k-1) \int_{\Omega} u^{k} c dx,
\end{align}
plugging the above formula into \er{star} we have
\begin{align}\label{starstar}
& \frac{d}{dt} \int_{\Omega}  u^k dx  + \frac{4(k-1)}{k}
\int_{\Omega} |\nabla u^{\frac{k}{2}} |^2 dx + k \int_{\Omega} u^\beta dx \int_{\Omega} u^{k+\alpha-1} dx+ (k-1) \int_{\Omega} u^{k} c dx  \nonumber\\
=& ~k \int_{\Omega} u^{k+\alpha-1} dx +(k-1) \int_{\Omega} u^{k+1} dx.
\end{align}
Omitting the last nonnegative term in the left hand side of \er{starstar} one has
\begin{align}\label{yuanshiguji}
& \frac{d}{dt} \int_{\Omega}  u^k dx  + \frac{4(k-1)}{k}
\int_{\Omega} |\nabla u^{\frac{k}{2}} |^2 dx + k \int_{\Omega} u^\beta dx \int_{\Omega} u^{k+\alpha-1} dx  \nonumber\\
\le & ~k \int_{\Omega} u^{k+\alpha-1} dx + (k-1) \int_{\Omega} u^{k+1} dx.
\end{align}
We will separate the proof into two cases $\alpha \ge 2$ and $\alpha<2$ in terms of the two nonnegative terms in the right hand side.
For $\alpha \ge 2$, $\int_{\Omega} u^{k+\alpha-1} dx$ is the dominant term and following procedures analogous to \cite{BCL15} we can obtain the boundedness of solutions in $L^k(1<k<\infty)$. While for the other case $\alpha<2$, the second term $\int_{\Omega} u^{k+1} dx$ dominates and we will use some appropriate Sobolev type of inequalities to show the boundedness of $L^k(1<k<\infty)$ norm.

\noindent {\bf Step 2 (A priori estimates for $\alpha \ge 2$).} For $\alpha \ge 2$, by Young's inequality one has
\begin{align}\label{HY}
 \int_{\Omega} u^{k+1} dx \le \int_{\Omega} u^{k+\alpha-1}dx+C(k,\alpha),
\end{align}
hence \er{yuanshiguji} becomes
\begin{align}\label{1yuanshiguji}
& \frac{d}{dt} \int_{\Omega}  u^k dx  + \frac{4(k-1)}{k}
\int_{\Omega} |\nabla u^{\frac{k}{2}} |^2 dx + k \int_{\Omega}u^\beta dx \int_{\Omega} u^{k+\alpha-1} dx  \nonumber\\
\le & ~C(k) \int_{\Omega} u^{k+\alpha-1} dx +C(k,\alpha).
\end{align}
Following similar procedures in \cite{BCL15} yields that if
\begin{align}\label{con1}
2 \le \alpha < 1+ 2\beta/n,
\end{align}
then it holds true that for any $1<k<\infty$ and any $0<t<T_{\max}$
\begin{align}
\|u(\cdot,t)\|_{L^k(\Omega)}^k \le C\left(\|u_0\|_{L^k(\Omega)}^k,k \right). \label{etaxiao}
\end{align}

\noindent {\bf Step 3 (A priori estimates for $\alpha<2$).} For $\alpha<2,$ similar to \er{HY}, \er{yuanshiguji} will be
\begin{align}\label{1yuanshiguji}
& \frac{d}{dt} \int_{\Omega}  u^k dx  + \frac{4(k-1)}{k}
\int_{\Omega} |\nabla u^{\frac{k}{2}} |^2 dx + k \int_{\Omega}u^\beta dx \int_{\Omega} u^{k+\alpha-1} dx  \nonumber\\
\le & ~C_2(k) \int_{\Omega} u^{k+1} dx +C(k,\alpha).
\end{align}
Letting $$~v=u^{k/2},q=\frac{2(k+1)}{k},r=\frac{2k'}{k}>1,C_0=\frac{2(k-1)}{k ~C_2(k) },C_1=\frac{1}{C_2(k)}$$
in Lemma \ref{lemmainterpolation} with
\begin{align}\label{kinterval}
k>\max \left\{ \frac{2}{p-2}, 1 \right\}
\end{align}
which is $q<p$ and
\begin{align}\label{kprime1}
k'>\frac{p}{p-2}
\end{align}
which is $\frac{q}{r}<\frac{2}{r}+1-\frac{2}{p},$ and $1<r<q<p$ equals to
\begin{align}\label{kprime}
\frac{k}{2}<k'<k+1<\frac{kp}{2},
\end{align}
one has
\begin{align}\label{12}
\int_{\Omega} u^{k+1} dx \le \frac{2(k-1)}{k ~C_2(k)} \| \nabla u^{\frac{k}{2}} \|_{L^2(\Omega)}^2+C(k) \|u\|_{L^{k'}(\Omega)}^{b}+\frac{1}{C_2(k)} \|u^{k/2} \|_{L^2(\Omega)}^2
\end{align}
with
$$
b=\frac{(1-\lambda)(k+1)}{1-\frac{\lambda (k+1)}{k}},~\lambda=\frac{\frac{k}{2k'}-\frac{k}{2(k+1)}}{\frac{k}{2k'}-\frac{1}{p}} \in(0,1).
$$
Combining \er{1yuanshiguji} with \er{12} we obtain that
\begin{align}\label{13}
& \frac{d}{dt}\int_{R^n} u^k dx+ k \int_{\Omega} u^\beta dx \int_{\Omega} u^{k+\alpha-1} dx+\frac{2(k-1)}{k} \|\nabla u^{\frac{k}{2}} \|_{L^2(\Omega)}^2 \nonumber \\
\le ~& C(k) \|u\|_{L^{k'}(\Omega)}^{b}+\|u\|_{L^k(\Omega)}^k+C(k,\alpha).
\end{align}
As
$$
\beta<k'<k+\alpha-1,
$$
we now use the following interpolation inequality
\begin{align}\label{14}
\|u\|_{L^{k'}(\Omega)}^{b} \le \left(\|u\|_{L^{k+\alpha-1}(\Omega)}^{k+\alpha-1} \|u\|_{L^\beta(\Omega)}^\beta \right)^{\frac{b \theta}{k+\alpha-1}} \|u\|_{L^\beta(\Omega)}^{b \left( 1-\theta-\frac{\theta \beta}{k+\alpha-1} \right)}
\end{align}
with
$$
\theta=\frac{\frac{1}{\beta}-\frac{1}{k'}}{\frac{1}{\beta}-\frac{1}{k+\alpha-1}} \in (0,1)
$$
to deal with \er{13}. Due to the arbitrariness of $k'$, we can take
$$
k'=\frac{k+\alpha-1+\beta}{2}
$$
such that
$$
1-\theta-\frac{\theta \beta}{k+\alpha-1}=0,
$$
and if
\begin{align}\label{17}
\frac{b \theta }{k+\alpha-1}<1,
\end{align}
then by using Young's inequality we can infer from \er{14} that
\begin{align}\label{23}
C(k) \|u\|_{L^{k'}(\Omega)}^{b} & \le C(k) \left(\|u\|_{L^{k+\alpha-1}(\Omega)}^{k+\alpha-1} \|u\|_{L^\beta(\Omega)}^\beta \right)^{\frac{b \theta}{k+\alpha-1}} \nonumber\\
& \le \frac{k}{4}\|u\|_{L^{k+\alpha-1}(\Omega)}^{k+\alpha-1} \|u\|_{L^\beta(\Omega)}^\beta +C(k,\alpha).
\end{align}
Now we discuss \er{17}. After a few computations, \er{17} is equivalent to
\begin{align}\label{19}
\frac{(2-\alpha)}{\beta}\left( \frac{k}{2}-\frac{k'}{p} \right)<(k+1-k')~\left( \frac{1}{2}-\frac{1}{p}-\frac{\alpha-1}{2 \beta} \right).
\end{align}
Introduce
\begin{align}
&A_0=\frac{1}{2}-\frac{1}{p}-\frac{\alpha-1}{2 \beta}>0, \\
&A_1=\frac{2-\alpha}{\beta}>0,
\end{align}
\er{19} can be written as
$$
A_0(k+1-k')>A_1 \left( \frac{k}{2}-\frac{k'}{p}  \right),
$$
that's
\begin{align}\label{20}
A_0+\left(A_0- \frac{A_1}{2} \right)k+\left( \frac{A_1}{p}-A_0 \right)k'>0.
\end{align}
If $A_0<\frac{A_1}{p}<\frac{A_1}{2}$, then \er{20} is
\begin{align}\label{61}
A_0+\left( \frac{A_1}{p}-A_0 \right)k'>\left(\frac{A_1}{2}-A_0\right)k.
\end{align}
Since $k'$ satisfies \er{kprime}, plugging $k'<\frac{kp}{2}$ into \er{61} follows
\begin{align}
\left(\frac{A_1}{2}-A_0\right)k< A_0+\left( \frac{A_1}{p}-A_0 \right)k'<A_0+\left( \frac{A_1}{p}-A_0 \right)\frac{kp}{2},
\end{align}
this is contrary to the fact \er{kinterval}. Otherwise if $A_0>\frac{A_1}{2}>\frac{A_1}{p},$ taking $k'>k/2$ and \er{20} into account one has
\begin{align}\label{247}
A_0+\left(A_0- \frac{A_1}{2} \right)k>\left( A_0-\frac{A_1}{p} \right)k'>\left( A_0-\frac{A_1}{p} \right) \frac{k}{2}.
\end{align}
Here we should remark that for $k'>\frac{p}{p-2}$ from \er{kprime1}, \er{247} reads
\begin{align}
A_0+\left(A_0- \frac{A_1}{2} \right)k>\left( A_0-\frac{A_1}{p} \right)k'>\left( A_0-\frac{A_1}{p} \right) \frac{p}{p-2},
\end{align}
this is just equivalent to \er{kinterval}, thus we only need to consider the case $k'>\frac{k}{2}$. Therefore, \er{17} holds true as long as $A_0-\frac{A_1}{2}>\frac{A_0}{2}-\frac{A_1}{2p},$ that's
\begin{align}\label{22}
\frac{(2-\alpha)}{\beta} \left(1-\frac{1}{p}\right) < \frac{1}{2}-\frac{1}{p}-\frac{\alpha-1}{2 \beta}.
\end{align}
On the other hand, letting
$$
v=u^{k/2}, q=2, 1\le r<2, C_0=\frac{k-1}{2k}, C_1=\frac{1}{2}
$$
in Lemma \ref{lemmainterpolation} obtains
$$
\int_{\Omega} u^k dx \le \frac{k-1}{2k}\|\nabla u^{k/2} \|_{L^2(\Omega)}^2+C(k,n)\|u\|_{L^{k_1}(\Omega)}^k+\frac{1}{2} \|u\|_{L^k(\Omega)}^k
$$
for $k_1=\frac{k r}{2}<k$ and thus
\begin{align}\label{xing}
\int_{\Omega} u^k dx \le \frac{k-1}{k}\|\nabla u^{k/2} \|_{L^2(\Omega)}^2+C(k,n)\|u\|_{L^{k_1}(\Omega)}^k.
\end{align}
Furthermore, for $\beta<k_1<k+\alpha-1$ we can take $k_1=\frac{\beta+\alpha-1+k}{2} \in (\beta,k)$ which is $k>\beta+\alpha-1$ such that
\begin{align}\label{xingxing}
\|u\|_{L^{\frac{\beta+\alpha-1+k}{2}}(\Omega)}^k \le \left( \|u\|_{L^{k+\alpha-1}(\Omega)}^{k+\alpha-1} \|u\|_{L^\beta(\Omega)}^\beta \right)^{\frac{k}{\beta+\alpha-1+k}}.
\end{align}
Taking \er{xing} and \er{xingxing} together and using Young's inequality yield
\begin{align}\label{xingxingxing}
\int_{\Omega} u^k dx \le \frac{k-1}{k}\|\nabla u^{k/2} \|_{L^2(\Omega)}^2+\frac{k}{4} \|u\|_{L^{k+\alpha-1}(\Omega)}^{k+\alpha-1} \|u\|_{L^\beta(\Omega)}^\beta +C(n,k).
\end{align}
Plugging \er{23} and \er{xingxingxing} into \er{13} one has that for $\beta+\alpha-1<k<\infty$
$$
\frac{d}{dt} \int_{\Omega} u^k dx +\int_{\Omega} u^k dx \le C(n,k)
$$
which follows
$$
\int_{\Omega} u^k dx \le C(\|u_0\|_{L^k(\Omega)},k).
$$
In addition, by virtue of Young's inequality, for any $1 \le k \le \beta+\alpha-1$
$$
\int_{\Omega}u^k dx \le \int_{\Omega} u^{\beta+\alpha}dx+C(n,k).
$$
Therefore we conclude that for all $1 \le k<\infty$
\begin{align}\label{etada}
\|u(\cdot,t)\|_{L^k(\Omega)} \le C\left(k,\|u_0\|_{L^k(\Omega)}\right).
\end{align}

\noindent {\bf Step 4 ($L^\infty$ estimates).} Based on the above arguments, firstly denoting $q_k=2^k+\beta+\alpha-1$ and taking $k=q_k$ in \er{yuanshiguji} we have
\begin{align}\label{qk}
& \frac{d}{dt} \int_{\Omega}  u^{q_k} dx  + \frac{4(q_k-1)}{q_k}
\int_{\Omega} |\nabla u^{\frac{q_k}{2}} |^2 dx + q_k \int_{\Omega}u^\beta dx \int_{\Omega} u^{q_k+\alpha-1} dx  \nonumber\\
\le & ~q_k \int_{\Omega} u^{q_k+\alpha-1} dx +(q_k-1) \int_{\Omega} u^{q_k+1} dx \nonumber \\
\le & ~q_k \int_{\Omega} u^{q_k+\alpha-1} dx +q_k \int_{\Omega} u^{q_k+1} dx.
\end{align}
Denoting
$$
s=\max\{\alpha,2\},
$$
by Young's inequality one has
\begin{align}
\int_{\Omega} u^{q_k+\alpha-1} dx \le \int_{\Omega} u^{q_k+s-1}+C(\alpha), \\
\int_{\Omega} u^{q_k+1} dx \le \int_{\Omega} u^{q_k+s-1}+C(\alpha).
\end{align}
Thus we infer from \er{qk} that
\begin{align}\label{qks}
& \frac{d}{dt} \int_{\Omega}  u^{q_k} dx  + \frac{4(q_k-1)}{q_k}
\int_{\Omega} |\nabla u^{\frac{q_k}{2}} |^2 dx + q_k \int_{\Omega} u^\beta dx \int_{\Omega} u^{q_k+\alpha-1} dx  \nonumber\\
\le & ~2q_k \int_{\Omega} u^{q_k+s-1} dx +q_k C(\alpha).
\end{align}
Now we apply Lemma \ref{lemmainterpolation} with
\begin{align}
v=u^{\frac{q_k}{2}},~~q=\frac{2(q_k+s-1)}{q_k},~~r=\frac{2 q_{k-1}}{q_k},~~C_0=\frac{1}{4 q_k},~~C_1=\frac{1}{4 q_k}~,
\end{align}
one has that
\begin{align}\label{qkinter}
\|u\|_{L^{q_k+s-1}(\Omega)}^{q_k+s-1} \le C(n) C_0^{\frac{-1}{\delta_1-1}} \left( \int_{\Omega} u^{q_{k-1}} dx \right)^{\gamma}+\frac{1}{4 q_k} \|\nabla u^{\frac{q_k}{2}} \|_{L^2(\Omega)}^2 + \frac{1}{4 q_k}\|u\|_{L^{q_k}(\Omega)}^{q_k}
\end{align}
where $\delta_1=\frac{q_k- 2 q_{k-1}/p}{q_k+s-1-q_{k-1}}=O(1)$, and from \er{con1} and \er{22} we can obtain that
\begin{align}
\gamma=1+\frac{q_k+s-1-q_{k-1}}{q_{k-1}-\frac{p(s-1)}{p-2}} \le 2.
\end{align}
Notice that $\frac{4(q_k -1)}{q_k} \ge 2$, substituting \er{qkinter} into \er{qks} follows
\begin{align}\label{qk1}
&\frac{d}{dt} \int_{\Omega} u^{q_k} dx +\frac{3}{2} \int_{\Omega} |\nabla u^{\frac{q_k}{2}}|^2 dx+ q_k \int_{\Omega}u^\beta dx \int_{\Omega} u^{q_k+\alpha-1}
dx \nonumber\\
\le & ~C(n) q_k^{\frac{\delta_1}{\delta_1 -1}} \left( \int_{\Omega} u^{q_{k-1}} dx \right)^{\gamma}+ \frac{1}{2} \|u\|_{L^{q_k}(\Omega)}^{q_k}+q_k C(\alpha).
\end{align}
On the other hand, using Lemma \ref{lemmainterpolation} with
\begin{align}
v=u^{\frac{q_k}{2}},~~q=2,~~r=\frac{2 q_{k-1}}{q_k},~~C_0=\frac{1}{4},~~C_1=\frac{1}{2}
\end{align}
we have
\begin{align}\label{qk2}
\frac{1}{2}\int_{\Omega} u^{q_k} dx  \le C(n) \left( \int_{\Omega} u^{q_{k-1}} dx \right)^{\frac{q_k}{q_{k-1}}} + \frac{1}{4}\int_{\Omega} |\nabla u^{\frac{q_k}{2}}|^2 dx.
\end{align}
In addition, using H\"{o}lder inequality with $\frac{1-\theta}{\beta}+\frac{\theta}{q_k+\alpha-1}=\frac{1}{q_{k-1}}$ and the fact
$$
q_{k-1}=\frac{q_k+\beta+\alpha-1}{2}
$$
we infer from Young's inequality that
\begin{align}\label{qk333}
C(n) \left( \int_{\Omega} u^{q_{k-1}} dx \right)^{\frac{q_k}{q_{k-1}}} & \le C(n) \left(\int_{\Omega} u^\beta dx \int_{\Omega} u^{q_k+\alpha-1} dx\right)^{\frac{q_k \theta }{q_k+\alpha-1}} \nonumber\\
& \le \int_{\Omega} u^\beta dx \int_{\Omega} u^{q_k+\alpha-1} dx +C(n,\alpha).
\end{align}
Combining \er{qk1}, \er{qk2} and \er{qk333} with the fact that $\gamma \le 2$ we have
\begin{align*}
\frac{d}{dt} \int_{\Omega} u^{q_k} dx+ \int_{\Omega} u^{q_k} dx  & \le C_1(n) q_k^{\frac{\delta_1}{\delta_1 -1}} \left( \int_{\Omega} u^{q_{k-1}} dx \right)^{\gamma}+C_2(n,\alpha)+q_k C_3(\alpha) \\
& \le C_1(n) q_k^{\frac{\delta_1}{\delta_1 -1}} \left( \int_{\Omega} u^{q_{k-1}} dx \right)^{\gamma}+ q_k^{\frac{\delta_1}{\delta_1 -1}} \left( C_2(n,\alpha)+C_3(\alpha) \right) \\
& =  C_1(n) q_k^{\frac{\delta_1}{\delta_1 -1}} \left( \int_{\Omega} u^{q_{k-1}} dx \right)^{\gamma}+ q_k^{\frac{\delta_1}{\delta_1 -1}} C(n,\alpha) \\
& \le \max\Big[ C_1(n),C(n,\alpha) \Big] q_k^{\frac{\delta_1}{\delta_1 -1}} \left[ \left( \int_{\Omega} u^{q_{k-1}} dx \right)^{\gamma}+1 \right] \\
& \le 2 \max \big[C_1(n),C(n,\alpha)\big] q_k^{\frac{\delta_1}{\delta_1 -1}} \max\left\{ \left( \int_{\Omega} u^{q_{k-1}} dx \right)^{2},1\right\}.
\end{align*}
Letting $r_0=\frac{\delta_1}{\delta_1 -1}$ and taking
\begin{align}
y_k(t)= \int_{\Omega} u^{q_k} dx,~\bar{a}=2 \max \big[C_1(n),C(n,\alpha)\big] (\alpha+\beta)^{r_0},~b=2
\end{align}
in Lemma \ref{odeinfinity} we obtain
\begin{align}\label{qk3}
\int_{\Omega} u^{q_k} dx \le (2 \bar{a})^{2^k-1} 2^{r_0 (2^{k+1}-k-2)} \max \left\{ \sup_{t\ge 0} \left(\int_{\Omega} u(t)^{q_0} dx\right)^{2^k}, K_0^{q_k} \right\}.
\end{align}
Here $K_0$ satisfies
\begin{align}
\int_{\Omega} u_0^{q_k} dx \le \left( \max \left\{ \|u_0\|_{L^{\beta+\alpha}(\Omega)},\|u_0\|_{L^\infty(\Omega)} \right\} \right)^{q_k} = K_0^{q_k}.
\end{align}
Taking the power $\frac{1}{q_k}$ to both sides of \er{qk3} and passing to the limit $k \to \infty$ one has
\begin{align}
\|u(t)\|_{L^\infty(\Omega)} \le 2 \bar{a} 2^{2 r_0} \max \left\{ \sup_{t\ge 0} \int_{\Omega} u(t)^{q_0} dx, K_0 \right\}.
\end{align}
Moreover, we infer from \er{etaxiao} and \er{etada} that
\begin{align}
\int_{\Omega} u(t)^{q_0} dx=\int_{\Omega} u(t)^{\beta+\alpha} dx \le C\left(\|u_0\|_{L^{\beta+\alpha}(\Omega)},\beta,\alpha\right) \le C\left( K_0,\beta,\alpha \right).
\end{align}
Consequently we obtain that for any $0<t<T_{\max}$
\begin{align}\label{infinityguji}
\|u(\cdot,t)\|_{L^\infty(\Omega)} \le C(K_0,\beta,\alpha).
\end{align}
Collecting \er{etaxiao},\er{etada} and \er{infinityguji} together we obtain the desired results. $\Box$ \\

\noindent {\bf Proof of Theorem \ref{thm0}:} Now we directly make use of the $L^\infty$ estimate \er{uniformbounded} and the blow-up criterion \er{blowup} to obtain that
\begin{align}
\|u(\cdot,t)\|_{L^\infty(\Omega)} \le C(\|u_0\|_{L^1(\Omega)}, \|u_0\|_{L^\infty(\Omega)})
\end{align}
for all $t \in (0,\infty).$ Thus completes the proof of Theorem \ref{thm0}. $\Box$

\end{document}